\documentclass[12pt,a4paper,twoside]{article}
\usepackage[english]{babel}
\usepackage{amssymb}
\usepackage{inputenc}
\newcommand{\ba}{\begin{array}}
\newcommand{\ea}{\end{array}}
\def\l {L^{0}}
\def\lt {L^{\infty}(\Omega)}
\usepackage{amssymb}
\begin{document}

\author{S. Albeverio $^{1},$ Sh. A. Ayupov $^{2},$ \ \ K. K.
Kudaybergenov  $^3$}
\title{\bf Derivations on the Algebra of Measurable Operators Affiliated with a   Type I von Neumann Algebra}

\maketitle

\begin{abstract}

Let $M$ be a type I von Neumann algebra with the center $Z,$  and
let $LS(M)$ be the algebra of all locally measurable operators
affiliated with $M.$ We prove that every $Z$-linear derivation on
$LS(M)$ is inner. In particular all $Z$-linear derivations on the
algebras of  measurable and respectively totally measurable
operators are spatial and implemented by  elements from $LS(M).$

\end{abstract}

\medskip
$^1$ Institut f\"{u}r Angewandte Mathematik, Universit\"{a}t Bonn,
Wegelerstr. 6, D-53115 Bonn (Germany); SFB 611, BiBoS; CERFIM
(Locarno); Acc. Arch. (USI),  e-mail: \emph{albeverio@uni-bonn.de}

$^2$ Institute of Mathematics and information  technologies,
Uzbekistan Academy of Science, F. Hodjaev str. 29, 100125,
Tashkent (Uzbekistan), e-mail: \emph{sh\_ayupov@mail.ru}

 $^{3}$ Institute of
Mathematics and information  technologies, Uzbekistan Academy of
Science, F. Hodjaev str. 29, 100125, Tashkent (Uzbekistan), e-mail:
\emph{karim2006@mail.ru}



\newpage

\section*{\center 1. Introduction}

The  present paper is devoted to  study   of derivations on the
algebra of locally measurable operators $LS(M)$ affiliated with a
type I  von Neumann algebra $M.$

Given a (complex) algebra $\mathcal{A},$ a linear operator  $D:
\mathcal{A}\rightarrow \mathcal{A}$
   is called a  \emph{derivation} if $D(xy)=D(x)y+xD(y)$ for all $x, y\in \mathcal{A}.$
   Each element $a\in \mathcal{A}$ generates a derivation $D_a:\mathcal{A}\rightarrow \mathcal{A}$
   defined as $D_a(x)=ax-xa,\,  x\in \mathcal{A}.$ Such derivations are called  \emph{inner} derivations.

It is well known that any derivation on a von Neumann algebra is
inner and therefore is norm continuous. But the properties of
derivations on the unbounded operator algebra $LS(M)$ seem to
be very far from being similar. Indeed, the results of
  \cite{Ber} and  \cite{Kus} show that in the commutative case when
 $M=L^{\infty}(\Omega, \Sigma, \mu),$ where $(\Omega, \Sigma, \mu)$
 is a non atomic measure space with a finite measure $\mu,$ the
 algebra $LS(M)\cong L^{0}(\Omega, \Sigma, \mu)$ of all classes of
   complex measurable functions on
$(\Omega, \Sigma, \mu)$ admits non zero derivations. It is clear
that these derivations are  discontinuous in the  measure topology
(i.~e. the  topology of convergence in measure),  and thus are non
inner. In order to avoid such pathological examples we consider
derivations on $LS(M)$ which are $Z$-linear, where $Z$ is the
center of the von Neumann algebra $M.$ The main result of the present
paper states that if $M$ is  a type I von Neumann algebra, then any
$Z$-linear derivation $D$ on the algebra $LS(M)$ is inner, i.~e.
$D_a(x)=ax-xa$ for an appropriate element $a\in LS(M).$

In Section 2 we give some preliminaries from the theory of
lattice-normed modules and Kaplansky~--- Hilbert modules over the
algebra of measurable functions and recall a result from
\cite{AyupKud} which gives a description of $\l$-linear
derivations on the algebra of all $\l$-bounded $\l$-linear
operators on a Banach~--- Kantorovich space over $\l=L^{0}(\Omega,
\Sigma, \mu).$

In Section 3 we prove that for any type I von Neumann algebra $M$
with the center $Z,$ every $Z$-linear derivation on the algebra
$LS(M)$ of locally measurable operators affiliated with $M$ is
inner. As a corollary we obtain that all $Z$-linear derivations on
the algebras of measurable and respectively totally measurable
operators affiliated with $M$ are spatial and implemented by
elements from $LS(M).$

\begin{center} {\bf 2. Preliminaries}
\end{center}

 Let   $(\Omega,\Sigma,\mu)$  be a measurable space and suppose
 that the measure $\mu$ has the  direct sum property, i.~e. there is a family $\{\Omega_{i}\}_{i\in
J}\subset\Sigma,\,0<\mu(\Omega_{i})<\infty,\,i\in J,$ such that for
any $A\in\Sigma,\,\mu(A)<\infty,$ there exist a countable subset
$J_{0 }\subset J$ and a set  $B$ with zero measure such that
$A=\bigcup\limits_{i\in J_{0}}(A\cap \Omega_{i})\cup B.$

 We
denote by  $\l=L^{0}(\Omega, \Sigma, \mu)$ the algebra of all
(classes of) complex measurable functions on $(\Omega, \Sigma, \mu)$
equipped with the topology of convergence in measure. Then $\l$ is a
complete metrizable commutative regular algebra with the unit
$\textbf{1}_\Omega$ given by $\textbf{1}_\Omega(\omega)=1, \,\omega\in\Omega.$

Recall that a net $\{\lambda_\alpha\}$ in $L^{0}$ $(o)$-converges
to $\lambda\in L^{0}$ if there exists a net $\{\xi_\alpha\}$
monotone decreasing to zero such that
$|\lambda_\alpha-\lambda|\leq\xi_\alpha$ for all $\alpha.$

 Denote
by $\nabla$ the complete Boolean algebra of all idempotents from
$\l,$ i. e. $\nabla=\{\tilde{\chi}_{A}: A\in\Sigma\},$ where
$\tilde{\chi}_{A}$ is  the element from $\l$ which contains the
characteristic function of the set $A.$

A complex linear space  $E$ is said to be normed  by $\l$ if there
is
 a map  $\|\cdot\|:E\longrightarrow \l$ such that for any  $x,y\in E, \lambda\in
\mathbb{C},$ the following conditions are fulfilled:

 1) $\|x\|\geq 0; \|x\|=0\Longleftrightarrow x=0$;

 2) $\|\lambda x\|=|\lambda|\|x\|$;

 3)  $\|x+y\|\leq\|x\|+\|y\|$.

 The  pair  $(E,\|\cdot\|)$ is called a lattice-normed  space over $\l.$
A lattice-normed  space  $E$ is called  $d$-decomposable, if for
any $x\in E$ with  $\|x\|=\lambda_{1}+\lambda_{2},$  $
\lambda_{i}\in\l,\,\lambda_i\geq0,\, i=1,2,\,
\lambda_{1}\lambda_{2}=0,$ there exist $x_{1}, x_{2}\in E$ such
that $x=x_{1}+x_{2}$ and $\|x_{i}\|=\lambda_{i},\,i=1,2$. A net
$(x_{\alpha})$ in $E$ is said to be $(bo)$-convergent to $x\in E$,
if the net $\{\|x_{\alpha}-x\|\}$ $(o)$-converges to zero in
 $\l.$
 A lattice-normed  space $E$ which is $d$-decomposable and  complete with respect to the  $(bo)$-convergence
  is called a
 \emph{Banach~--- Kantorovich space}.

 It is known  that every Banach~--- Kantorovich space $E$ over
 $\l$ is a module over $\l$ and  $\|\lambda x\|=|\lambda|\|x\|$
 for all $\lambda\in\l,\, x\in E$ (see \cite{Kusr}).

Any Banach~--- Kantorovich space  $E$ over $\l$ is orthocomplete, i.~e.
given any net  $(x_{\alpha})\subset E$ and a partition of the unit
$(\pi_{\alpha})$ in $\nabla$ the series
$\sum\limits_{\alpha}\pi_{\alpha}x_{\alpha}$ $(bo)$-converges in
$E.$

 A module  $F$ over  $\l$ is said to be  finite-generated, if there
 exist elements
        $x_{1},x_{2},...,x_{n}$ in $F$ such that  any  $x\in F$ can be decomposed as
        $x=\lambda_{1}x_{1}+...+\lambda_{n}x_{n}$ for appropriate  $\lambda_{i}\in
       \l\, (i=\overline{1,n}).$
       The elements $x_{1},x_{2},...,x_{n}$ are called  generators of $F.$
       We denote by $d(F)$ the minimal number of  generators of $F.$
       A  module  $F$ over  $\l$ is  called $\sigma$-finite-generated,
       if there exists
       a partition $(\pi_{\alpha})_{\alpha\in A}$ of the unit in $\nabla$ such that
       $\pi_{\alpha}F$ is finite-generated for any $\alpha.$
         A finite-generated module  $F$ over  $\l$ is called homogeneous of type  $n$, if $d(eF)=n$ for every
          nonzero $e\in\nabla.$

Let  $\mathcal{K}$ be a module over $\l$. A map $\langle
\cdot,\cdot\rangle:\mathcal{K}\times \mathcal{K}\rightarrow\l$ is
called  an $\l$-valued inner product, if for all $x,y,z\in
\mathcal{K},\,\lambda\in\l,$ it satisfies the following
conditions:

1) $\langle x,x\rangle\geq0$; $\langle x,x\rangle=0\Leftrightarrow
x=0$;

2) $\langle x,y\rangle=\overline{\langle y,x\rangle}$;

3) $\langle \lambda x,y\rangle=\lambda\langle x,y\rangle$;

4) $\langle x+y,z\rangle=\langle x,z\rangle+\langle y,z\rangle$.

If $\langle \cdot,\cdot\rangle:\mathcal{K}\times
\mathcal{K}\rightarrow\l$ is an $\l$-valued inner product, then $
\|x\|=\sqrt{\langle x,x \rangle} $ defines  an $\l$-valued norm on
$\mathcal{K}.$ The  pair $(\mathcal{K},\langle \cdot,\cdot\rangle)$
is called a \emph{Kaplansky~--- Hilbert module} over $\l,$ if
$(\mathcal{K},\|\cdot\|)$ is a Banach~--- Kantorovich space over  $\l$
(see \cite{Kusr}).

Let  $X$ be a Banach space.  A map  $s:\Omega\rightarrow X$ is said
to be simple, if
$s(\omega)=\sum\limits_{k=1}^{n}\chi_{A_{k}}(\omega)c_k,$ where
$A_k\in\Sigma, A_i\cap A_j=\emptyset, \,i\neq j,\, \,c_k\in
X,\,k=\overline{1, n},\, n\in\mathbb{N}.$ A map $u:\Omega\rightarrow
X$ is said to be measurable, if there is a sequence  $(s_n)$ of
simple maps such that $\|s_n(\omega)-u(\omega)\|\rightarrow0$
almost everywhere on any $A\in\sum$ with $\mu(A)<\infty.$

Let $\mathcal{L}(\Omega, X)$ be the set of all measurable maps
from $\Omega$ into $X,$ and let $L^{0}(\Omega, X)$ denote the
space of all equivalence classes  in $\mathcal{L}(\Omega, X)$ with
respect to the equality almost everywhere. Denote by $\hat{u}$ the
equivalence class from $L^{0}(\Omega, X)$ which contains the
measurable map $u\in \mathcal{L}(\Omega, X).$ Further we shall
identify the element $u\in \mathcal{L}(\Omega, X)$ and the class
$\hat{u}.$ Note that the function  $\omega \rightarrow
\|u(\omega)\|$
       is measurable for any $u\in \mathcal{L}(\Omega, X).$ The equivalence class containing the function
              $\|u(\omega)\|$ is denoted by
       $\|\hat{u}\|$. For  $\hat{u}, \hat{v}\in L^{0}(\Omega, X), \lambda\in\l$ put
$\hat{u}+\hat{v}=\widehat{u(\omega)+v(\omega)},
\lambda\hat{u}=\widehat{\lambda(\omega) u(\omega)}.$

It is known  \cite{Kusr} that $(L^{0}(\Omega, X), \|\cdot\|)$ is a
Banach~--- Kantorovich space over $\l.$

Put $$L^{\infty}(\Omega)=\{f\in L^{0}: \exists c\in \mathbb{R}, c>0,
|f|\leq c\textbf{1}\}$$ and $$L^{\infty}(\Omega, X)=\{x\in
L^{0}(\Omega, X):\|x\|\in L^{\infty}(\Omega)\}.$$  Then
$L^{\infty}(\Omega, X)$ is a Banach space with respect to  the norm
$$\|x\|_{\infty}=\|\|x\|\|_{L^{\infty}(\Omega)},\,x\in L^{\infty}(\Omega, X).$$

If $H$ is a Hilbert space, then $L^{0}(\Omega, H)$ can be equipped
with an  $\l$-valued inner product $\langle x,
y\rangle=\widehat{(x(\omega), y(\omega))},$ where $(\cdot, \cdot)$
is the inner product on $H.$ Then $(L^{0}(\Omega, H),\langle
\cdot,\cdot\rangle)$ is a Kaplansky~--- Hilbert module over $\l$
and  $(L^{\infty}(\Omega, H),\langle \cdot,\cdot\rangle)$ is a
Kaplansky~--- Hilbert module over $L^{\infty}(\Omega).$

Let  $E$ be a  Banach~--- Kantorovich space over $\l.$
 An operator   $T: E\rightarrow E$
 is   $\l$-linear if  $T(\lambda_1 x_1 +\lambda_2 x_2)=\lambda_1 T(x_1)+\lambda_2
 T(x_2)$
for all  $\lambda_1, \lambda_2\in \l ,x_1, x_2\in E.$ An $\l$-linear
operator
 $T:E\rightarrow E$ is said to be  $\l$-bounded if there exists an element $c\in\l$
 such that $\|T(x)\|\leq
 c\|x\|$ for any $x\in E.$  For an  $\l$-bounded $\l$-linear
operator $T$  we put $\|T\|=\sup\{\|T(x)\|:\|x\|\leq \textbf{1} \}.$

An $\l$-bounded $\l$-linear operator $T:E\rightarrow E$ is called
finite-generated (respectively $\sigma$-finite-generated) if
$T(E)=\{T(x):x\in E\}$ is a finite-generated (respectively
$\sigma$-finite-generated) submodule in $E.$

Denote by  $B(E)$   the  algebra of all
 $\l$-boun\-ded $\l$-linear operators
on  $E$
 and  let $F_{\sigma}(E)$  be the set of all
$\sigma$-finite-generated  operators on $E.$

Let
 $B(L^{\infty}(\Omega, H))$ be
 the set of all  $L^{\infty}(\Omega)$-bounded $L^{\infty}(\Omega)$-linear operators
 on $L^{\infty}(\Omega, H).$

Put
  $$B(L^{0}(\Omega, H))_b=\{x\in B(L^{0}(\Omega, H)): \|x\|\in L^{\infty}(\Omega)\}.$$
  Note that   the correspondence
    $$x\mapsto x|_{L^{\infty}(\Omega, H)}$$
    gives  a $\ast$-isomorphism between the $\ast$-algebras
    $B(L^{0}(\Omega, H))_b$ and $B(L^{\infty}(\Omega, H)).$
 Further we shall identify  $B(L^{0}(\Omega, H))_b$ with
    $B(L^{\infty}(\Omega, H))$ (i.~e.  the operator $x$ from $B(L^{0}(\Omega,
    H))_b$ is identified
     with  its restriction $x|_{L^{\infty}(\Omega, H)}$).

     We shall conclude this section with the following theorem
     from \cite{AyupKud}, which is necessary for the proof of the
     main
result of the present paper.

 \textbf{Theorem 2.1.} \cite{AyupKud}.  \emph{Let  $E$ be a  Banach~--- Kantorovich space over $\l$
 and  let $D:B(E)\rightarrow B(E)$ be an
$\l$-linear derivation.
 Then there is  $T\in B(E)$ such that $$D(A)=TA-AT$$ for all} $A\in B(E).$

 \begin{center} {\bf 3. Derivations on  the algebra of locally measurable operators for type I von Neumann algebras}
\end{center}

Let  $B(H)$ be the algebra of all bounded linear operators on a
Hilbert space  $H$ and let  $M$ be a von Neumann algebra in $B(H).$
 Denote by
$\mathcal{P}(M)$ the lattice of projections in  $M.$

A linear subspace  $\mathcal{D}$ in  $H$ is said to be affiliated
with  $M$ (denoted as  $\mathcal{D}\eta M$), if
$u(\mathcal{D})\subset \mathcal{D}$ for any unitary operator $u$
from the commutant
$$M'=\{y\in B(H):xy=yx, \,\forall x\in M\}$$ of the algebra $M.$

A linear operator  $x$ on  $H$ with the domain  $\mathcal{D}(x)$
is said to be affiliated with  $M$ (denoted as  $x\eta M$) if
$u(\mathcal{D}(x))\subset \mathcal{D}(x)$ and $ux(\xi)=xu(\xi)$
for every unitary operator $u\in M'$ and all  $\xi\in
\mathcal{D}(x).$

A linear subspace $\mathcal{D}$ in $H$ is said to be strongly
dense in  $H$ with respect to the von Neumann algebra  $M,$ if

1) $\mathcal{D}\eta M;$

2) there exists a sequence of projections
$\{p_n\}_{n=1}^{\infty}\subset P(M),$  such that
$p_n\uparrow\textbf{1},$ $p_n(H)\subset \mathcal{D}$ and
$p^{\perp}_n=\textbf{1}-p_n$ is finite in  $M$ for all
$n\in\mathbb{N},$ where $\textbf{1}$ is the identity in $M.$

A closed linear operator  $x$ on a Hilbert space $H$ is said to be
\emph{measurable} with respect to the von Neumann algebra  $M,$ if
 $x\eta M$ and $\mathcal{D}(x)$ is strongly dense in  $H.$ Denote by
 $S(M)$ the set of all measurable operators affiliated with
 $M$ (see \cite{Seg}).

A closed linear operator $x$  on a Hilbert space  $H$  is said to
be \emph{locally measurable} with respect to the von Neumann
algebra $M,$ if $x\eta M$ and there exists a sequence
$\{z_n\}_{n=1}^{\infty}$ of central projections in $M$ such that
$z_n\uparrow\textbf{1}$ and $z_nx \in S(M)$ for all
$n\in\mathbb{N}$  (see \cite{Yea1}).

It  is known \cite{Yea1} that the set $LS(M)$ of all locally
measurable operators affiliated with  $M$ forms a unital
$\ast$-algebra with respect to the strong algebraic operations and
the natural involution. Moreover  $S(M)$ is a solid
$\ast$-subalgebra in $LS(M).$

The following result describes one of the most important
properties of the algebra  $LS(M)$ (see \cite{ChilLit},
\cite{Sai}).

\textbf{Proposition 3.1.} \emph{Suppose that the von Neumann
algebra  $M$ is the  $C^{\ast}$-product of the von Neumann
algebras  $M_i,$ $i\in I,$ where $I$ is an arbitrary set of
indices, i~.e.}
$$M=\sum\limits_{i\in I}^{\oplus}M_i=
\{\{x_i\}_{i\in I}:x_i\in M_i, i\in I, \sup\limits_{i\in
I}\|x_i\|_{M_i}<\infty\}$$ \emph{with coordinate-wise algebraic
operations and involution and with the $C^{\ast}$-norm
$\|\{x_i\}_{i\in I}\|_{ M}=\sup\limits_{i\in I}\|x_i\|_{M_i}.$
Then the algebra  $LS(M)$ is $\ast$-isomorphic to the algebra
$\prod\limits_{i\in I}LS(M_i)$ (with the coordinate-wise
operations and involution).}

This proposition implies that given any family  $\{z_i\}_{i\in I}$
of mutually orthogonal central projections in $M$ with
$\bigvee\limits_{i\in I}z_i=\textbf{1}$ and any family
$\{x_i\}\subset LS(M)$ there exists a unique element $x\in LS(M)$
such that $z_i x=z_i x_i$ for all $i\in I.$ This element is
denoted by  $x=\sum\limits_{i\in I}z_i x_i.$

Recall that a von Neumann algebra $M$ is an algebra of \emph{type I}
if it is isomorphic to a von Neumann algebra with an abelian
commutant.

The main result of the present paper is the following

\textbf{Theorem 3.2.} \emph{Let   $M$ be a type I von Neumann
algebra with the center $Z.$ Then every  $Z$-linear derivation on
the algebra  $LS(M)$ is inner.}

The main tool in the proof of this theorem is the decomposition of
the given von Neumann algebra in to  the direct sum of homogeneous
components and the representation of homogeneous type I von
Neumann algebras as algebras of bounded $Z$-linear operators on a
Kaplansky~--- Hilbert module over the center $Z$ of the given von
Neumann algebra.

For details we refer the reader to the monograph of A.G. Kusraev
\cite{Kusr}.

First, let us consider the case of a homogeneous type I von
Neumann algebra.

Let $M$ be a homogeneous von Neumann algebra of type  I$_\alpha,$
where  $\alpha$ is a cardinal number, and let $L^{\infty}(\Omega)$
be the center of  $M.$ Then it is known that  $M$ is
$\ast$-isomorphic to the algebra  $B(L^{\infty}(\Omega, H)),$
where $\dim H=\alpha$ (see \cite{Kap},  \cite{Kusr}).

Let  $(\Omega, \Sigma, \mu)$  be a measure space such that the
measure $\mu$ has the direct sum property, and let
$\{\Omega_i\}_{i\in\mathbb{N}}$ be a countable partition of the
set  $\Omega$ into measurable subsets. Let  $H_{n_{i}}$ denote a
finite dimensional Hilbert space with the dimension  $n_i, i\in
\mathbb{N}.$

Put
$$\sum\limits_{i\in\mathbb{N}}^{\oplus}L^{\infty}(\Omega_{i}, H_{n_i})=
\{\{\varphi_i\}_{i\in\mathbb{N}}: \varphi_i\in L^{\infty}(\Omega_{i}, H_{n_i}),
\{\|\varphi_i\|_{i}\}_{i\in\mathbb{N}}\in
 \sum\limits_{i\in\mathbb{N}}^{\oplus}L^{\infty}(\Omega_{i})\},$$ where
$\|\cdot\|_{i}$ is
$L^{\infty}(\Omega_{i})$-valued norm on  $L^{\infty}(\Omega_{i},
H_{n_i}).$  Equipped with coordinate-wise algebraic
 operations and inner product the set
$\sum\limits_{i\in\mathbb{N}}^{\oplus}L^{\infty}(\Omega_{i},
H_{n_i})$ becomes a Kaplansky~--- Hilbert module over
$\sum\limits_{i\in\mathbb{N}}^{\oplus}L^{\infty}(\Omega_{i})\cong
L^{\infty}(\Omega)$ ($\cong$ standing for $\ast$-isomorphism).

Similarly, the set
 $$\prod\limits_{i\in\mathbb{N}}L^{0}(\Omega_{i},
H_{n_i})$$ is a Kaplansky~--- Hilbert module over
$\prod\limits_{i\in\mathbb{N}}L^{0}(\Omega_{i})\cong \l.$

Note that  $\prod\limits_{i\in\mathbb{N}}L^{0}(\Omega_{i},
H_{n_i})$ and
$\sum\limits_{i\in\mathbb{N}}^{\oplus}L^{\infty}(\Omega_{i},
H_{n_i})$ are the general forms of $\sigma$-finite-generated
Kaplansky~--- Hilbert modules over $L^{0}$ and over
$L^{\infty}(\Omega),$ respectively.

  It should be noted that each finite von Neumann algebra of type
   I with the center  $L^{\infty}(\Omega)$ is  $\ast$-isomorphic to an appropriate
   algebra of the form
$B(\sum\limits_{i\in\mathbb{N}}^{\oplus}L^{\infty}(\Omega_{i},
H_{n_i})),$ where  $\dim H_{n_i}=n_i, n_i\in\mathbb{N},$   and
moreover
$$LS(M)=S(M)\cong
B(\prod\limits_{i\in\mathbb{N}}(L^{0}(\Omega_{i}, H_{n_i}))).$$

\textbf{Proposition 3.3.} \emph{A projection  $p\in
B(L^{\infty}(\Omega, H))$ is finite if and only if  $p$ is
$\sigma$-finite-generated.}

Proof. Let  $p\in B(L^{\infty}(\Omega, H))$ be a projection. Since
 $pB(L^{\infty}(\Omega, H))p\cong B(p(L^{\infty}(\Omega, H))),$
 $p$ is finite if and only if the algebra
$B(p(L^{\infty}(\Omega, H)))$ is finite. But this exactly means
that  $p(L^{\infty}(\Omega, H))$ is $\sigma$-finite-generated
module, i.~e.   $p$ is  $\sigma$-finite-generated projection. The
proof is complete.

Denote by  $B(L^{\infty}(\Omega, H))+F_{\sigma}(L^{0}(\Omega, H))$
the $\ast$-subalgebra in  $B(L^{0}(\Omega, H))$ which consists of
elements of the form  $x=x_1+x_2,$ with  $x_1\in
B(L^{\infty}(\Omega, H)),$ $x_2\in F_{\sigma}(L^{0}(\Omega, H)).$

\textbf{Proposition 3.4.} \emph{The algebra
$S(B(L^{\infty}(\Omega, H)))$ is
 $\ast$-isomorphic to the  algebra $B(L^{\infty}(\Omega, H))+F_{\sigma}(L^{0}(\Omega, H)).$}

Proof. Take    $x\in S(B(L^{\infty}(\Omega, H)))$ and let
 $\{e_{\lambda}\}_{\lambda\in\mathbb{R}}$ be the spectral family
 of the element  $|x|.$ By   \cite[Propsition
 1]{ChilLit} there exists  $\lambda_0$ such that
$e_{\lambda_{0}}^{\perp}$ is a finite projection.  We have
$$xe_{\lambda_{0}}^{\perp}=e_{\lambda_{0}}^{\perp}xe_{\lambda_{0}}^{\perp}\in
e_{\lambda_{0}}^{\perp} S(B(L^{\infty}(\Omega,
H)))e_{\lambda_{0}}^{\perp}$$ and
$$e_{\lambda_{0}}^{\perp} S(B(L^{\infty}(\Omega, H)))e_{\lambda_{0}}^{\perp}\cong
S(e_{\lambda_{0}}^{\perp} B(L^{\infty}(\Omega,
H))e_{\lambda_{0}}^{\perp})\cong
S(B(e_{\lambda_{0}}^{\perp}(L^{\infty}(\Omega, H))).$$ By
Proposition  3.3  $e_{\lambda_{0}}^{\perp}$ is a
$\sigma$-finite-generated projection and therefore
$$S(B(e_{\lambda_{0}}^{\perp}(L^{\infty}(\Omega, H)))\cong
B(e_{\lambda_{0}}^{\perp}(L^{0}(\Omega, H))).$$ Thus, under the
obtained  $\ast$-isomorphism
$$e_{\lambda_{0}}^{\perp} S(B(L^{\infty}(\Omega,
H)))e_{\lambda_{0}}^{\perp} \cong
B(e_{\lambda_{0}}^{\perp}(L^{0}(\Omega, H)))$$ the element  $x
e_{\lambda_{0}}^{\perp}$ corresponds to some
  $\sigma$-finite-generated operator from $B(e_{\lambda_{0}}^{\perp}(L^{0}(\Omega, H))),$
  which is
  denoted by  $\widetilde{x e_{\lambda_{0}}^{\perp}}.$ Since
  $x e_{\lambda_{0}}\in B(L^{\infty}(\Omega, H)),$ we have that the mapping
 $$\Phi: x\mapsto x e_{\lambda_{0}}+\widetilde{xe_{\lambda_{0}}^{\perp}}$$
 gives a $\ast$-embedding of the algebra  $S(B(L^{\infty}(\Omega, H)))$
 into $B(L^{\infty}(\Omega, H))+F_{\sigma}(L^{0}(\Omega, H)).$

Now let  $x\in F_{\sigma}(L^{0}(\Omega, H)).$ Take a
$\sigma$-finite-generated projection
 $p\in B(L^{\infty}(\Omega, H))$ such that $x=pxp.$
 Then $x=pxp\in pB(L^{0}(\Omega, H))p.$

 Since $p$ is a
$\sigma$-finite-generated projection,  $B(p(L^{\infty}(\Omega,
 H)))$ is a finite von Neumann  algebra. Hence
 $$pB(L^{0}(\Omega, H))p\cong B(p(L^{0}(\Omega,
 H)))\cong S(B(p(L^{\infty}(\Omega, H)))\cong$$
 $$\cong S(pB(L^{\infty}(\Omega, H))p)\cong pS(B(L^{\infty}(\Omega, H)))p.$$
 Thus $pB(L^{0}(\Omega, H))p\cong
pS(B(L^{\infty}(\Omega, H)))p.$ Hence, under this
$\ast$-isomorphism the operator $x$ corresponds to some element
from $pS(B(L^{\infty}(\Omega, H)))p\subset S(B(L^{\infty}(\Omega,
H))),$ and therefore the mapping $\Phi$ is surjective.

 This means that   $\Phi$ is a $\ast$-isomorphism between $S(B(L^{\infty}(\Omega, H)))$
  and  $B(L^{\infty}(\Omega, H))+F_{\sigma}(L^{0}(\Omega, H)).$ The
proof is complete.

\textbf{Proposition   3.5.} \emph{The algebras
$LS(B(L^{\infty}(\Omega, H)))$ and  $B(L^{0}(\Omega, H))$ are
$\ast$-isomorphic.}

Proof. Let us show that the  $\ast$-isomorphism  $\Phi$ between
$S(B(L^{\infty}(\Omega, H)))$ and  $B(L^{\infty}(\Omega,
H))+F_{\sigma}(L^{0}(\Omega, H))$  can be extended to a
$\ast$-isomorphism between  $LS(B(L^{\infty}(\Omega, H)))$
  and   $B(L^{0}(\Omega, H)).$

  Let $x\in LS(B(L^{\infty}(\Omega,
H))).$  Consider a sequence  $\{z_n\}$ of central projections such
that  $z_n\uparrow \textbf{1}$ and $xz_n\in S(B(L^{\infty}(\Omega,
H)))$ for all $n\in\mathbb{N}.$ Put $\pi_1=z_1,$ $\pi_n=z_n\wedge
z_{n-1}^{\perp}$ $n\geq2,$ and
$$\Psi(x)=(bo)-\sum\limits\limits_{n\in\mathbb{N}}\pi_n\Phi(\pi_n x_n).$$
Since  $B(L^{0}(\Omega, H))$ is orthocomplete,  $\Psi$ is an
imbedding of the algebra  $LS(B(L^{\infty}(\Omega, H)))$ into the
algebra $B(L^{0}(\Omega, H)).$

Let now $y\in B(L^{0}(\Omega, H)).$ Take a sequence  $\{z_n\}$ of
mutually orthogonal central projections such that
$\bigvee\limits_{n\in\mathbb{N}}z_n =\textbf{1}$ and $\|y\| z_n\in
L^{\infty}(\Omega).$ Then  $y z_n\in B(L^{\infty}(\Omega, H)).$
Put $x=\sum\limits_{n\in \mathbb{N}}z_n y.$  Then $\Psi(x)=y$ and
therefore  $\Psi$ is a surjective map. This means that  $\Psi$ is
a $\ast$-isomorphism between  $LS(B(L^{\infty}(\Omega, H)))$
  and  $B(L^{0}(\Omega, H)).$ The proof is complete.

  The assertions of Propositions  3.4 and 3.5 become more clear in the following
  case of homogeneous type  I$_\infty$ von Neumann algebra with the discrete center.

Let $M$ be the $C^{\ast}$-product of a countable family of copies
of the von Neumann algebra  $B(H)$ with  $\dim H=\infty,$
 i.~e.
 $$M=\sum\limits_{n\in\mathbb{N}}^{\oplus}B(H).$$
Proposition 3.5 and 3.4 imply that
$$LS(M)\cong\prod\limits_{n\in\mathbb{N}}B(H)$$
and
$$S(M)\cong M+\prod\limits_{n\in\mathbb{N}}F(H),$$
where  $F(H)$ is the ideal of finite-dimensional operators from
$B(H).$

Now let us consider general von Neumann algebras of type I.

It is well-known \cite{Tak} that if $M$ is a type I von Neumann
algebra
 then there is a unique (cardinal-indexed) orthogonal
family of projections $(q_{\alpha})_{\alpha\in
I}\subset\mathcal{P}(M)$ with $\sum\limits_{\alpha\in
I}q_{\alpha}=\textbf{1}$ such that  $q_{\alpha}M$ is a homogeneous
type I$_\alpha$ von Neumann algebra,
 i.~e. $q_\alpha M\cong B(L^{\infty}(\Omega_{\alpha}, H_{\alpha})),\, \dim H_\alpha=\alpha,$ and
$$M\cong\sum\limits_{\alpha\in I}^{\oplus}B(L^{\infty}(\Omega_{\alpha}, H_{\alpha})).$$
Note that if  $L^{\infty}(\Omega)$ is the center of the von
Neumann algebra  $M$ then
 $q_\alpha L^{\infty}(\Omega)\cong L^{\infty}(\Omega_{\alpha})$
 for all $\alpha\in I.$

The product
$$\prod\limits_{\alpha\in I}L^{0}(\Omega_{\alpha}, H_{\alpha})$$
is a Kaplansky~--- Hilbert module over $L^{0}$ with respect to the
coordinate-wise algebraic operations and inner product.

 The product
$$\prod\limits_{\alpha\in I}B(L^{0}(\Omega_{\alpha},
H_{\alpha}))$$ with the coordinate-wise algebraic operations and
involution forms a  $\ast$-algebra and moreover
$$\prod\limits_{\alpha\in I}B(L^{0}(\Omega_{\alpha},
H_{\alpha}))\cong B(\prod\limits_{\alpha\in
I}L^{0}(\Omega_{\alpha}, H_{\alpha})).\eqno (1)$$

Now Propositions  3.1,  3.5 and the isomorphism  (1) imply

\textbf{Proposition 3.6.} \emph{For $M\cong\sum\limits_{\alpha\in
I}^{\oplus}B(L^{\infty}(\Omega_{\alpha}, H_{\alpha}))$ the algebra
$LS(M)$ is
 $\ast$-isomorphic to the algebra  $B(\prod\limits_{\alpha\in I}L^{0}(\Omega_{\alpha},
H_{\alpha})).$}

From Proposition  3.6 it follows that if $M$ is a type I von
Neumann algebra then given any $x\in LS(M)$ there exists a
sequence $\{z_n\}$ of mutually orthogonal central projections with
$\bigvee\limits_{n\in\mathbb{N}}z_n =\textbf{1}$ and $z_n x\in M$
for all $n\in\mathbb{N};$

Proof of the Theorem 3.2. Let $D:LS(M)\rightarrow LS(M)$ be a
$Z$-linear derivation and $\lambda\in\l.$ Take a sequence $(e_n)$
of projections in $Z$ such that $e_n\lambda\in \lt=Z$ and
$e_n\uparrow\textbf{1}.$ Then for any $n$ and $x\in LS(M)$ we have
   $e_n D(\lambda x)=D(e_n \lambda x)=e_n\lambda D(x)$ and therefore
    $D(\lambda x)=\lambda D(x),$ i. e. any $Z$-linear derivation on $LS(M)$ is also $\l$-linear.
  By Proposition 3.6 and Theorem 2.1 we have $D$ is inner. The proof is complete.

\textbf{Remark 3.7.} The condition on the derivation to be
$Z$-linear is crucial in general. This follows from the examples
of non zero derivations on the commutative algebra $L^{0}(0;
1)\cong LS(L^{\infty}(0; 1))$ given in \cite{Ber} (see also \cite{Kus}).
These derivations are not inner, and moreover they are not
continuous in the measure topology. Another example of
discontinuous (and hence non inner) derivation is the following
non commutative generalization of the above one.

\textbf{Example 3.8.} Let $\delta$ be any of   non zero
derivations on $L^{0}(0; 1)$ constructed in \cite{Ber}. Consider
the von Neumann algebra $M=L^{\infty}(0; 1)\overline{\otimes}
M_n(\mathbb{C}),$ which can be identified with the algebra of all
$n\times n$ matrices $(f_{i, j})_{i, j=1}^{n}$ with entries from
$L^{\infty}(0; 1).$ Then the algebra $LS(M)$ is $\ast$-isomorphic
to  the algebra $M_n(L^{0}(0; 1))$ of all $n\times n$ matrices
with entries $f_{i, j}$ from the algebra $L^{0}(0; 1).$

Define the mapping $D_\delta:M_n(L^{0}(0; 1)) \rightarrow M_n(L^{0}(0; 1))$ by
$$D_\delta ((f_{i, j})_{i, j=1}^{n})=(\delta(f_{i, j}))_{i, j=1}^{n}.$$
Then it is easy to check that $D_\delta$ is a derivation on
$M_n(L^{0}(0; 1))$
 which is  not $Z$-linear (where
$Z=L^{\infty}(0; 1)),$  and that $D_\delta$ is discontinuous and hence can
not be inner.

Now let  $M$ be a type  I von Neumann algebra and let
$\mathcal{A}$ be an arbitrary subalgebra of $LS(M)$ containing
$M.$ Consider a derivation  $D:\mathcal{A}\rightarrow\mathcal{A}$
and let us show that  $D$ can be extended to a derivation
$\tilde{D}$ on the whole $LS(M).$

For an arbitrary element  $x\in LS(M)$ take a sequence  $\{z_n\}$
of mutually orthogonal central projections with
$\bigvee\limits_{n\in\mathbb{N}}z_n =\textbf{1}$ and  $z_n x\in M$
for all $n\in\mathbb{N}.$

Put $$\tilde{D}(x)=\sum z_n D(z_n x).\eqno (2)$$

Since every derivation  $D:\mathcal{A}\rightarrow\mathcal{A}$ is
identically zero on central projections of  $M,$ the equality (2)
gives a well-defined $Z$-linear derivation
$\tilde{D}:LS(M)\rightarrow LS(M)$ which coincides with $D$ on
$\mathcal{A}.$ By  Theorem 3.2 the derivation  $\tilde{D}$ is
inner and therefore  $D$ is a spatial derivation on $\mathcal{A},$
i.~e. there exists an element $a\in LS(M)$ such that
$$D(x)=ax-xa$$
for all $x\in \mathcal{A}.$

Therefore we obtain the following

\textbf{Theorem 3.9.} \emph{Let     $M$ be a type  I von Neumann
algebra with the center $Z,$ and let   $\mathcal{A}$ be an
arbitrary subalgebra in $LS(M)$ containing  $M.$ Then any
$Z$-linear derivation on  $\mathcal{A}$  is spatial and
implemented by an element of  $LS(M).$}

Now let  $\tau$ be a faithful normal semi-finite trace on the von
Neumann algebra  $M.$ Recall that a closed linear operator  $x$ is
said to be  $\tau$-measurable (or totally measurable) with respect
to the von Neumann algebra  $M,$ if $x\eta M$ and its domain
$\mathcal{D}(x)$ is  $\tau$-dense in $H$ (i.~e.
$\mathcal{D}(x)\eta M$ and given any  $\varepsilon>0$ there exists
a projection $p\in\mathcal{P}(M)$ such that  $p(H)\subset
\mathcal{D}(x)$ and $\tau(p^{\perp})\leq\varepsilon$).

The set  $S(M, \tau)$ of all $\tau$-measurable operators with
respect to $M$ is a solid $\ast$-subalgebra in $S(M)$ (see
\cite{ChilLit}). Therefore Theorem  3.9 implies

\textbf{Corollary  3.10.} \emph{Let    $M$ be a type  I von
Neumann algebra with the center  $Z$ and let  $D$ be a $Z$-linear
derivation on  $S(M)$ or on  $S(M, \tau).$ Then  $D$ is spatial
and implemented by an element from  $LS(M).$}

Now let  $M$ be a type I   von Neumann algebra with the atomic
center $Z,$ and let $\{q_i\}_{i\in I}$ be the set of all atoms
from $Z.$ Consider a derivation $D$ on $S(M, \tau).$ Since $q_i
Z\cong q_i \mathbb{C}$ for all $i\in I,$ we have $q_i D(\lambda
x)=D(q_i \lambda x)=q_i \lambda D(x)$ for all $i\in I,\,
\lambda\in Z,\,x\in S(M, \tau).$ Thus $D(\lambda x)=\lambda D(x)$
for any $\lambda\in Z.$ This means that in the case of atomic $Z$
any derivation on $S(M, \tau)$ is automatically $Z$-linear. From
this and from  Corollary 3.10 we have the following result, which
is a strengthening of a result of Weigt \cite{Wei}.

\textbf{Corollary 3.11.} \emph{If  $M$ is a von Neumann algebra
with the atomic lattice of projections, then every   derivation on
the algebra $S(M, \tau)$ is spatial, and in particular is
continuous in the measure topology.} \vspace{1cm}

\textbf{Acknowledgments.} \emph{The second and third named authors
would like to acknowledge the hospitality of the $\,$ "Institut
f\"{u}r Angewandte Mathematik",$\,$ Universit\"{a}t Bonn (Germany).
This work is supported in part by the DFG 436 USB 113/10/0-1 project
(Germany) and the Fundamental Research  Foundation of the Uzbekistan
Academy of Sciences.}

\emph{The authors are indebted to the reviewer for useful
comments.}

\end{document}